# EMPIRICAL STUDY IN FINITE CORRELATION COEFFICIENT IN TWO PHASE ESTIMATION

M. Khoshnevisan, Lecturer, Griffith University, School of Accounting and Finance, Australia.
F. Kaymarm, Assistant Professor, Massachusetts Institute of Technology, Department of Mechanical Engineering, USA; currently at Sharif University, Tehran, Iran.
H. P. Singh, R Singh, Professors of Statistics, Vikram University, Department of Mathematics and Statistics, India.
F. Smarandache, Associate Professor, Department of Mathematics, University of New Mexico, Gallup, USA.

## ABSTRACT

This paper proposes a class of estimators for population correlation coefficient when information about the population mean and population variance of one of the variables is not avaliable but information about these parameters of another variable (auxiliary) is avaliable, in two phase sampling and analyzes its properties. Optimum estimator in the class is identified with its variance formula. The estimators of the class involve unknown constants whose optimum values depend on unknown population parameters.Following Singh (1982) and Srivastava and Jhajj (1983), it has been shown that when these population parameters are replaced by their consistent estimates the resulting class of estimators has the same asymptotic variance as that of optimum estimator. An empirical study is carried out to demonstrate the performance of the constructed estimators.

## 1. Introduction

Consider a finite population U= {1,2,..,i,..N}. Let y and x be the study and auxiliary variables taking values $y_i$ and $x_i$ respectively for the ith unit. The correlation coefficient between y and x is defined by

$$\rho_{yx} = S_{yx} / (S_y S_x) \qquad (1.1)$$

where

$$S_{yx} = (N-1)^{-1} \sum_{i=1}^{N} (y_i - \overline{Y})(x_i - \overline{X}),\ S_x^2 = (N-1)^{-1} \sum_{i=1}^{N} (x_i - \overline{X})^2,\ S_y^2 = (N-1)^{-1} \sum_{i=1}^{N} (y_i - \overline{Y})^2,$$

$$\overline{X} = N^{-1} \sum_{i=1}^{N} x_i,\ \overline{Y} = N^{-1} \sum_{i=1}^{N} y_i .$$

Based on a simple random sample of size *n* drawn without replacement,

$(x_i, y_i)$, $i = 1,2,\ldots,n$; the usual estimator of $\rho_{yx}$ is the corresponding sample correlation coefficient :

$$r = s_{yx}/(s_x s_y) \qquad (1.2)$$

where $s_{yx} = (n-1)^{-1} \sum_{i=1}^{n}(y_i - \bar{y})(x_i - \bar{x})$, $s_x^2 = (n-1)^{-1} \sum_{i=1}^{n}(x_i - \bar{x})^2$

$s_y^2 = (n-1)^{-1} \sum_{i=1}^{n}(y_i - \bar{y})^2$, $\bar{y} = n^{-1} \sum_{i=1}^{n} y_i$, $\bar{x} = n^{-1} \sum_{i=1}^{n} x_i$.

The problem of estimating $\rho_{yx}$ has been earlier taken up by various authors including Koop (1970), Gupta et. al. (1978, 79), Wakimoto (1971), Gupta and Singh (1989), Rana (1989) and Singh et. al. (1996) in different situations. Srivastava and Jhajj (1986) have further considered the problem of estimating $\rho_{yx}$ in the situations where the information on auxiliary variable *x* for all units in the population is available. In such situations, they have suggested a class of estimators for $\rho_{yx}$ which utilizes the known values of the population mean $\bar{X}$ and the population variance $S_x^2$ of the auxiliary variable *x*.

In this paper, using two – phase sampling mechanism, a class of estimators for $\rho_{yx}$ in the presence of the available knowledge ($\bar{Z}$ and $S_z^2$) on second auxiliary variable z is considered, when the population mean $\bar{X}$ and population variance $S_x^2$ of the main auxiliary variable *x* are not known.

## 2. The Suggested Class of Estimators

In many situations of practical importance, it may happen that no information is available on the population mean $\bar{X}$ and population variance $S_x^2$, we seek to estimate the population correlation coefficient $\rho_{yx}$ from a sample 's' obtained through a two-phase selection. Allowing simple random sampling without replacement scheme in each phase, the two- phase sampling scheme will be as follows:

(i) The first phase sample $s^*$ $(s^* \subset U)$ of fixed size $n_1$, is drawn to observe only *x* in order to furnish a good estimates of $\bar{X}$ and $S_x^2$.

(ii) Given $s^*$, the second- phase sample s $(s \subset s^*)$ of fixed size *n* is drawn to observe *y* only.

Let

$\bar{x} = (1/n) \sum_{i \in s} x_i$, $\bar{y} = (1/n) \sum_{i \in s} y_i$, $\bar{x}^* = (1/n_1) \sum_{i \in s^*} x_i$, $s_x^2 = (n-1)^{-1} \sum_{i \in s}(x_i - \bar{x})^2$,

$s_x^{*2} = (n_1 - 1)^{-1} \sum_{i \in s^*}(x_i - \bar{x}^*)^2$.

We write $u = \bar{x}/\bar{x}^*$, $v = s_x^2 / s_x^{*2}$. Whatever be the sample chosen let $(u,v)$ assume values in a bounded closed convex subset, R, of the two-dimensional real space containing the point (1,1). Let *h* $(u, v)$ be a function of *u* and *v* such that

$$h(1,1) = 1 \qquad (2.1)$$

and such that it satisfies the following conditions:
1. The function *h* $(u,v)$ is continuous and bounded in R.
2. The first and second partial derivatives of *h(u,v)* exist and are continuous and bounded in R.

Now one may consider the class of estimators of $\rho_{yx}$ defined by
$$\hat{\rho}_{hd} = r\, h(u,v) \qquad (2.2)$$
which is double sampling version of the class of estimators

$$\tilde{r}_t = r\, f(u^*, v^*)$$

Suggested by Srivastava and Jhajj (1986), where $u^* = \bar{x}/\bar{X}$, $v^* = s_x^2/S_x^2$ and $(\bar{X}, S_x^2)$ are known.

Sometimes even if the population mean $\bar{X}$ and population variance $S_x^2$ of x are not known, information on a cheaply ascertainable variable z, closely related to x but compared to x remotely related to y, is available on all units of the population. This type of situation has been briefly discussed by, among others, chand (1975), Kiregyera (1980, 84).

Following Chand (1975) one may define a chain ratio- type estimator for $\rho_{yx}$ as

$$\hat{\rho}_{1d} = r\left(\frac{\bar{x}^*}{\bar{x}}\right)\left(\frac{\bar{Z}}{\bar{z}^*}\right)\left(\frac{s_x^{*2}}{s_x^2}\right)\left(\frac{S_z^2}{s_z^{*2}}\right) \qquad (2.3)$$

where the population mean $\bar{Z}$ and population variance $S_z^2$ of second auxiliary variable z are known, and

$$\bar{z}^* = (1/n_1)\sum_{i \in s^*} z_i, \quad s_z^{*2} = (n_1 - 1)^{-1}\sum_{i \in s^*}(z_i - \bar{z}^*)^2$$

are the sample mean and sample variance of z based on preliminary large sample $s^*$ of size $n_1\ (>n)$.

The estimator $\hat{\rho}_{1d}$ in (2.3) may be generalized as

$$\hat{\rho}_{2d} = r\left(\frac{\bar{x}}{\bar{x}^*}\right)^{\alpha_1}\left(\frac{s_x^2}{s_x^{*2}}\right)^{\alpha_2}\left(\frac{\bar{z}^*}{\bar{Z}}\right)^{\alpha_3}\left(\frac{s_z^{*2}}{S_z^2}\right)^{\alpha_4} \qquad (2.4)$$

where $\alpha_i's$ (i=1,2,3,4) are suitably chosen constants.

Many other generalization of $\hat{\rho}_{1d}$ is possible. We have, therefore, considered a more general class of $\rho_{yx}$ from which a number of estimators can be generated.

The proposed generalized estimators for population correlation coefficient $\rho_{yx}$ is defined by
$$\hat{\rho}_{td} = r\, t(u,v,w,a) \qquad (2.5)$$
where $w = \bar{z}^*/\bar{Z}$, $a = s_z^{*2}/S_z^2$ and *t(u,v,w,a)* is a function of *(u,v,w,a)* such that
$$t\,(1,1,1,1)=1 \qquad (2.6)$$
Satisfying the following conditions:
(i)     Whatever be the samples ($s^*$ and s) chosen, let *(u,v,w,a)* assume values in a closed convex subset S, of the four dimensional real space containing the point P=(1,1,1,1).
(ii)    In S, the function *t(u,v,w,a)* is continuous and bounded.
(iii)   The first and second order partial derivatives of *t(u,v,w, a)* exist and are continuous and bounded in S
To find the bias and variance of $\hat{\rho}_{td}$ we write

$s_y^2 = S_y^2(1+e_1), \bar{x} = \bar{X}(1+e_1), \bar{x}^* = \bar{X}(1+e_1^*), s_x^2 = S_x^2(1+e_2)$

$s_x^{*2} = S_x^2(1+e_2^*), \bar{z}^* = \bar{Z}(1+e_3^*), s_z^{*2} = S_z^2(1+e_4^*), s_{yx} = S_{yx}(1+e_s)$

such that $E(e_0) = E(e_1) = E(e_2) = E(e_5) = 0$ and $E(e_i^*) = 0 \ \forall \ i = 1,2,3,4$, and ignoring the finite population correction terms, we write to the first degree of approximation

$E(e_0^2) = (\delta_{400} - 1)/n, E(e_1^2) = C_x^2/n, E(e_1^{*2}) = C_x^2/n_1, E(e_2^2) = (\delta_{040} - 1)/n$,

$E(e_2^{*2}) = (\delta_{040} - 1)/n_1, E(e_3^{*2}) = C_z^2/n_1, E(e_4^{*2}) = (\delta_{004} - 1)/n_1$,

$E(e_5^2) = \{(\delta_{220}/\rho_{yx}^2) - 1\}/n, E(e_0 e_1) = \delta_{210} C_x/n, E(e_0 e_1^*) = \delta_{210} C_x/n_1$,

$E(e_0 e_2) = (\delta_{220} - 1)/n, E(e_0 e_2^*) = (\delta_{220} - 1)/n_1, E(e_0 e_3^*) = \delta_{201} C_z/n_1$,

$E(e_0 e_4^*) = (\delta_{202} - 1)/n_1, E(e_0 e_5) = \{(\delta_{310}/\rho_{yx}) - 1\}/n$,

$E(e_1 e_1^*) = C_x^2/n_1, E(e_1 e_2) = \delta_{030} C_x/n, E(e_1 e_2^*) = \delta_{030} C_x/n_1$,

$E(e_1 e_3^*) = \rho_{xz} C_x C_z/n_1, E(e_1 e_4^*) = \delta_{012} C_x/n_1, E(e_1 e_5) = (\delta_{120} C_x/\rho_{yx})/n$,

$E(e_1^* e_2) = \delta_{030} C_x/n_1, E(e_1^* e_2^*) = \delta_{030} C_x/n_1, E(e_1^* e_3^*) = \rho_{xz} C_x C_z/n_1$,

$E(e_1^* e_4^*) = \delta_{012} C_x/n_1, E(e_1^* e_5) = (\delta_{120} C_x/\rho_{yx})/n_1$,

$E(e_2 e_2^*) = (\delta_{040} - 1)/n_1, E(e_2 e_3^*) = \delta_{021} C_z/n_1, E(e_2 e_4^*) = (\delta_{022} - 1)/n_1$,

$E(e_2 e_5) = \{(\delta_{130}/\rho_{yx}) - 1\}/n, \ E(e_2^* e_3^*) = \delta_{021} C_z/n_1$,

$E(e_2^* e_4^*) = (\delta_{022} - 1)/n_1, E(e_2^* e_5) = \{(\delta_{130}/\rho_{yx}) - 1\}/n_1$,

$E(e_3^* e_4^*) = \delta_{003} C_z/n_1, E(e_3^* e_5) = (\delta_{111} C_z/\rho_{yx})/n_1$,

$E(e_4^* e_5) = \{(\delta_{112}/\rho_{yx}) - 1\}/n_1$.

where

$\delta_{pqm} = \mu_{pqm}/(\mu_{200}^{p/2} \mu_{020}^{q/2} \mu_{002}^{m/2}), \ \mu_{pqm} = (1/N)\sum_{i=1}^{N}(y_i - \bar{Y})^p (x_i - \bar{X})^q (z_i - \bar{Z})^m$, $(p,q,m)$ being non-negative integers.

To find the expectation and variance of $\hat{\rho}_{td}$, we expand $t(u,v,w,a)$ about the point $P = (1,1,1,1)$ in a second-order Taylor's series, express this value and the value of r in terms of e's. Expanding in powers of e's and retaining terms up to second power, we have

$$E(\hat{\rho}_{td}) = \rho_{yx} + o(n^{-1}) \qquad (2.7)$$

which shows that the bias of $\hat{\rho}_{td}$ is of the order $n^{-1}$ and so up to order $n^{-1}$, mean square error and the variance of $\hat{\rho}_{td}$ are same.

Expanding $(\hat{\rho}_{td} - \rho_{yx})^2$, retaining terms up to second power in e's, taking expectation and using the above expected values, we obtain the variance of $\hat{\rho}_{td}$ to the first degree of approximation, as

$$Var(\hat{\rho}_{td}) = Var(r) + (\rho_{yx}^2/n)[C_x^2 t_1^2(P) + (\delta_{040} - 1)t_2^2(P) - At_1(P) - Bt_2(P) + 2\delta_{030}C_x t_1(P)t_2(P)]$$
$$- (\rho_{yx}^2/n_1)[C_x^2 t_1^2(P) + (\delta_{040} - 1)t_2^2(P) - C_z^2 t_3^2(P) - (\delta_{004} - 1)t_4^2(P) - At_1(P) -$$
$$Bt_2(P) + Dt_3(P) + Ft_4(P) + 2\delta_{030}C_x t_1(P)t_2(P) - 2\delta_{003}C_z t_3(P)t_4(P)]$$

(2.8)

where $t_1(P)$, $t_2(P)$, $t_3(P)$ and $t_4(P)$ respectively denote the first partial derivatives of $t(u,v,w,a)$ white respect to $u,v,w$ and $a$ respectively at the point $P=(1,1,1,1)$,

$$Var(r) = (\rho_{yx}^2/n)[\delta_{220}/\rho_{yx}^2 + (1/4)(\delta_{040} + \delta_{400} + 2\delta_{220}) - \{(\delta_{130} + \delta_{310})/\rho_{yx}\}] \quad (2.9)$$

$$A = \{\delta_{210} + \delta_{030} - 2(\delta_{120}/\rho_{yx})\}C_x, B = \{\delta_{220} + \delta_{040} - 2(\delta_{130}/\rho_{yx})\},$$
$$D = \{\delta_{201} + \delta_{021} - 2(\delta_{111}/\rho_{yx})\}C_z, F = \{\delta_{202} + \delta_{022} - 2(\delta_{112}/\rho_{yx})\}$$

Any parametric function $t(u,v,w,a)$ satisfying (2.6) and the conditions (1) and (2) can generate an estimator of the class (2.5).

The variance of $\hat{\rho}_{td}$ at (2.6) is minimized for

$$\left. \begin{array}{l} t_1(P) = \dfrac{[A(\delta_{040} - 1) - B\delta_{030}C_x]}{2C_x^2(\delta_{040} - \delta_{030}^2 - 1)} = \alpha \text{ (say)}, \\[2mm] t_2(P) = \dfrac{(BC_x^2 - A\delta_{030}C_x)}{2C_x^2(\delta_{040} - \delta_{030}^2 - 1)} = \beta \text{ (say)}, \\[2mm] t_3(P) = \dfrac{[D(\delta_{004} - 1) - F\delta_{003}C_z]}{2C_z^2(\delta_{004} - \delta_{030}^2 - 1)} = \gamma \text{ (say)}, \\[2mm] t_4(P) = \dfrac{(C_z^2 F - D\delta_{003}C_z)}{2C_z^2(\delta_{004} - \delta_{003}^2 - 1)} = \delta \text{ (say)}, \end{array} \right\} \quad (2.10)$$

Thus the resulting (minimum) variance of $\hat{\rho}_{td}$ is given by

$$\min. Var(\hat{\rho}_{td}) = Var(r) - (\frac{1}{n} - \frac{1}{n_1})\rho_{yx}^2 [\frac{A^2}{4C_x^2} + \frac{\{(A/C_x)\delta_{030} - B\}^2}{4(\delta_{040} - \delta_{030}^2 - 1)}]$$
$$- (\rho_{yx}^2/n_1)\left[\frac{D^2}{4C_z^2} + \frac{\{(D/C_z)\delta_{003} - F\}^2}{4(\delta_{004} - \delta_{003}^2 - 1)}\right]$$

(2.11)

It is observed from (2.11) that if optimum values of the parameters given by (2.10) are used, the variance of the estimator $\hat{\rho}_{td}$ is always less than that of r as the last two terms on the right hand sides of (2.11) are non-negative.

Two simple functions $t(u,v,w,a)$ satisfying the required conditions are

$t(u,v,w,a) = 1 + \alpha_1(u - 1) + \alpha_2(v - 1) + \alpha_3(w - 1) + \alpha_4(a - 1)$

$t(u,v,w,a) = u^{\alpha_1} v^{\alpha_2} w^{\alpha_3} a^{\alpha_4}$

and for both these functions $t_1(P) = \alpha_1$, $t_2(P) = \alpha_2$, $t_3(P) = \alpha_3$ and $t_4(P) = \alpha_4$. Thus one should use optimum values of $\alpha_1, \alpha_2, \alpha_3$ and $\alpha_4$ in $\hat{\rho}_{td}$ to get the minimum variance. It is to be noted that the estimated $\hat{\rho}_{td}$ attained the minimum variance only when the optimum values of the constants $\alpha_i$ (i=1,2,3,4), which are functions of unknown population parameters, are known. To use such estimators in practice, one has to use some guessed values of population parameters obtained either through past experience or through a pilot sample survey. It may be further noted that even if the values of the constants used in the estimator are not exactly equal to their optimum values as given by (2.8) but are close enough, the resulting estimator will be better than the conventional estimator, as has been illustrated by Das and Tripathi (1978, Sec.3).

If no information on second auxiliary variable z is used, then the estimator $\hat{\rho}_{td}$ reduces to $\hat{\rho}_{hd}$ defined in (2.2). Taking $z \equiv 1$ in (2.8), we get the variance of $\hat{\rho}_{hd}$ to the first degree of approximation, as

$$Var(\hat{\rho}_{hd}) = Var(r) + \left(\frac{1}{n} - \frac{1}{n_1}\right)\rho_{yx}^2 \left[C_x^2 h_1^2(1,1) + (\delta_{040} - 1)h_2^2(1,1) - Ah_1(1,1) - Bh_2(1,1) + 2\delta_{030} C_x h_1(1,1)h_2(1,1)\right]$$

(2.12)

which is minimized for

$$h_1(1,1) = \frac{[A(\delta_{040} - 1) - B\delta_{030} C_x]}{2C_x^2(\delta_{040} - \delta_{030}^2 - 1)}, \quad h_2(1,1) = \frac{(BC_x^2 - A\delta_{030} C_x)}{2C_x^2(\delta_{040} - \delta_{030}^2 - 1)} \quad (2.13)$$

Thus the minimum variance of $\hat{\rho}_{hd}$ is given by

$$\min. Var(\hat{\rho}_{hd}) = Var(r) - \left(\frac{1}{n} - \frac{1}{n_1}\right)\rho_{yx}^2 \left[\frac{A^2}{4C_x^2} + \frac{\{(A/C_x)\delta_{030} - B\}^2}{4(\delta_{040} - \delta_{030}^2 - 1)}\right] \quad (2.14)$$

It follows from (2.11) and (2.14) that

$$\min. Var(\hat{\rho}_{td}) - \min. Var(\hat{\rho}_{hd}) = \left(\rho_{yx}^2/n_1\right)\left[\frac{D^2}{4C_z^2} + \frac{\{(D/C_z)\delta_{003} - F\}^2}{4(\delta_{004} - \delta_{003}^2 - 1)}\right] \quad (2.15)$$

which is always positive. Thus the proposed estimator $\hat{\rho}_{td}$ is always better than $\hat{\rho}_{hd}$.

## 3. A Wider Class of Estimators

In this section we consider a class of estimators of $\rho_{yx}$ wider than (2.5) given by

$$\hat{\rho}_{gd} = g(r,u,v,w,a) \quad (3.1)$$

where $g(r,u,v,w,a)$ is a function of $r,u,v,w,a$ and such that

$$g(\rho,1,1,1,1) = \rho \text{ and } \left[\frac{\partial g(\cdot)}{\partial r}\right]_{(\rho,1,1,1)} = 1$$

Proceeding as in section 2, it can easily be shown, to the first order of approximation, that the minimum variance of $\hat{\rho}_{gd}$ is same as that of $\hat{\rho}_{td}$ given in (2.11).

It is to be noted that the difference-type estimator
$r_d = r + \alpha_1 (u-1) + \alpha_2 (v-1) + \alpha_3 (w-1) + \alpha_4 (a-1)$, is a particular case of $\hat{\rho}_{gd}$, but it is not the member of $\hat{\rho}_{td}$ in (2.5).

## 4. Optimum Values and Their Estimates

The optimum values $t_1(P) = \alpha$, $t_2(P) = \beta$, $t_3(P) = \gamma$ and $t_4(P) = \delta$ given at (2.10) involves unknown population parameters. When these optimum values are substituted in (2.5), it no longer remains an estimator since it involves unknown ($\alpha, \beta, \gamma, \delta$), which are functions of unknown population parameters, say, $\delta_{pqm}$ (p, q, m= 0,1,2,3,4), $C_x$, $C_z$ and $\rho_{yx}$ itself. Hence it is advisable to replace them by their consistent estimates from sample values. Let $(\hat{\alpha}, \hat{\beta}, \hat{\gamma}, \hat{\delta})$ be consistent estimators of $t_1(P), t_2(P), t_3(P)$ and $t_4(P)$ respectively, where

$$\hat{t}_1(P) = \hat{\alpha} = \frac{[\hat{A}(\hat{\delta}_{040} - 1) - \hat{B}\hat{\delta}_{030}\hat{C}_x]}{2\hat{C}_x^2(\hat{\delta}_{040} - \hat{\delta}_{030}^2 - 1)}, \qquad \hat{t}_2(P) = \hat{\beta} = \frac{[\hat{B}\hat{C}_x^2 - \hat{A}\hat{\delta}_{030}\hat{C}_x]}{2\hat{C}_x^2(\hat{\delta}_{040} - \hat{\delta}_{030}^2 - 1)},$$

$$\hat{t}_3(P) = \hat{\gamma} = \frac{[\hat{D}(\hat{\delta}_{004} - 1) - \hat{F}\hat{\delta}_{003}\hat{C}_z]}{2\hat{C}_z^2(\hat{\delta}_{004} - \hat{\delta}_{003}^2 - 1)}, \qquad \hat{t}_4(P) = \hat{\delta} = \frac{[\hat{C}_z^2\hat{F} - \hat{D}\hat{\delta}_{003}\hat{C}_z]}{2\hat{C}_z^2(\hat{\delta}_{004} - \hat{\delta}_{003}^2 - 1)},$$

(4.1)

with

$$\hat{A} = [\hat{\delta}_{210} + \hat{\delta}_{030} - 2(\hat{\delta}_{120}/r)]\hat{C}_x, \qquad \hat{B} = [\hat{\delta}_{220} + \hat{\delta}_{040} - 2(\hat{\delta}_{130}/r)],$$

$$\hat{D} = [\hat{\delta}_{201} + \hat{\delta}_{021} - 2(\hat{\delta}_{111}/r)]\hat{C}_z, \qquad \hat{F} = [\hat{\delta}_{202} + \hat{\delta}_{022} - 2(\hat{\delta}_{112}/r)],$$

$$\hat{C}_x = s_x/\bar{x}, \quad \hat{C}_z = s_z/\bar{z}, \quad \hat{\delta}_{pqm} = \hat{\mu}_{pqm}/\left(\hat{\mu}_{200}^{p/2} \hat{\mu}_{020}^{q/2} \hat{\mu}_{002}^{m/2}\right)$$

$$\hat{\mu}_{pqm} = (1/n)\sum_{i=1}^{n}(y_i - \bar{y})^p(x_i - \bar{x})^q(z_i - \bar{z})^m$$

$$\bar{z} = (1/n)\sum_{i=1}^{n} z_i, \quad s_x^2 = (n-1)^{-1}\sum_{i=1}^{n}(x_i - \bar{x})^2, \quad \bar{x} = (1/n)\sum_{i=1}^{n} x_i,$$

$$r = s_{yx}/(s_y s_x), \quad s_y^2 = (n-1)^{-1}\sum_{i=1}^{n}(y_i - \bar{y})^2, \quad s_z^2 = (n-1)^{-1}\sum_{i=1}^{n}(x_i - \bar{z})^2.$$

We then replace $(\alpha, \beta, \gamma, \delta)$ by $(\hat{\alpha}, \hat{\beta}, \hat{\gamma}, \hat{\delta})$ in the optimum $\hat{\rho}_{td}$ resulting in the estimator $\hat{\rho}_{td}^*$ say, which is defined by

$$\hat{\rho}_{td}^* = rt^*(u, v, w, a, \hat{\alpha}, \hat{\beta}, \hat{\gamma}, \hat{\delta}), \qquad (4.2)$$

where the function $t^*(U)$, $U = (u, v, w, a, \hat{\alpha}, \hat{\beta}, \hat{\gamma}, \hat{\delta})$ is derived from the the function $t(u,v,w,a)$ given at (2.5) by replacing the unknown constants involved in it by the consistent estimates of optimum values. The condition (2.6) will then imply that

$$t^*(P^*) = 1 \qquad (4.3)$$

where $\qquad P^* = (1,1,1,1, \alpha, \beta, \gamma, \delta)$

We further assume that

$$t_1(P^*) = \left.\frac{\partial t^*(U)}{\partial u}\right|_{U=P^*} = \alpha, \qquad t_2(P^*) = \left.\frac{\partial t^*(U)}{\partial v}\right|_{U=P^*} = \beta$$

$$t_3(P^*) = \frac{\partial t^*(U)}{\partial w}\bigg]_{U=P^*} = \gamma, \qquad t_4(P^*) = \frac{\partial t^*(U)}{\partial a}\bigg]_{U=P^*} = \delta \qquad (4.4)$$

$$t_5(P^*) = \frac{\partial t^*(U)}{\partial \hat{\alpha}}\bigg]_{U=P^*} = 0 \qquad t_6(P^*) = \frac{\partial t^*(U)}{\partial \hat{\beta}}\bigg]_{U=P^*} = 0$$

$$t_7(P^*) = \frac{\partial t^*(U)}{\partial \hat{\gamma}}\bigg]_{U=P^*} = 0 \qquad t_8(P^*) = \frac{\partial t^*(U)}{\partial \hat{\delta}}\bigg]_{U=P^*} = 0$$

Expanding $t^*(U)$ about $P^* = (1,1,1,1,\alpha,\beta,\gamma,\delta)$, in Taylor's series, we have

$$\hat{\rho}_{td}^* = r[t^*(P^*) + (u-1)t_1^*(P^*) + (v-1)t_2^*(P^*) + (w-1)t_3^*(P^*) + (a-1)t_4^*(P^*) + (\hat{\alpha}-\alpha)t_5^*(P^*)$$
$$+ (\hat{\beta}-\beta)t_6^*(P^*) + (\hat{\gamma}-\gamma)t_7^*(P^*) + (\hat{\delta}-\delta)t_8^*(P^*) + \text{second order terms}]$$
(4.5)

Using (4.4) in (4.5) we have

$$\hat{\rho}_{td}^* = r[1 + (u-1)\alpha + (v-1)\beta + (w-1)\gamma + (a-1)\delta + \text{second order terms}] \qquad (4.6)$$

Expressing (4.6) in term of e's squaring and retaining terms of e's up to second degree, we have

$$(\hat{\rho}_{td}^* - \rho_{yx})^2 = \rho_{yx}^2 [\frac{1}{2}(2e_5 - e_0 - e_2) + \alpha(e_1 - e_1^*) + \beta(e_2 - e_2^*) + \gamma\, e_3^* + \delta\, e_4^*]^2 \qquad (4.7)$$

Taking expectation of both sides in (4.7), we get the variance of $\hat{\rho}_{td}^*$ to the first degree of approximation, as

$$Var(\hat{\rho}_{td}^*) = Var(r) - (\frac{1}{n} - \frac{1}{n_1})\rho_{yx}^2 \left[\frac{A^2}{4C_x^2} + \frac{\{(A/C_x)\delta_{030} - B\}^2}{4(\delta_{040} - \delta_{030}^2 - 1)}\right]$$
$$+ (\rho_{yx}^2/n_1)\left[\frac{D^2}{4C_z^2} + \frac{\{(D/C_z)\delta_{003} - F\}^2}{4(\delta_{004} - \delta_{003}^2 - 1)}\right]$$
(4.8)

which is same as (2.11), we thus have established the following result.

**Result 4.1:** If optimum values of constants in (2.10) are replaced by their consistent estimators and conditions (4.3) and (4.4) hold good, the resulting estimator $\hat{\rho}_{td}^*$ has the same variance to the first degree of approximation, as that of optimum $\hat{\rho}_{td}$.

**Remark 4.1:** It may be easily examined that some special cases:

(i) $\hat{\rho}^*_{td1} = r u^{\hat{\alpha}} v^{\hat{\beta}} w^{\hat{\gamma}} a^{\hat{\delta}}$, (ii) $\hat{\rho}^*_{td2} = r \dfrac{\{1+\hat{\alpha}(u-1)+\hat{\gamma}(w-1)\}}{\{1-\hat{\beta}(v-1)-\hat{\delta}(a-1)\}}$

(iii) $\hat{\rho}^*_{td3} = r[1+\hat{\alpha}(u-1)+\hat{\beta}(u-1)+\hat{\gamma}(w-1)+\hat{\delta}(a-1)]$

(iv) $\hat{\rho}^*_{td4} = r[1-\hat{\alpha}(u-1)-\hat{\beta}(u-1)-\hat{\gamma}(w-1)-\hat{\delta}(a-1)]^{-1}$

of $\hat{\rho}^*_{td}$ satisfy the conditions (4.3) and (4.4) and attain the variance (4.8).

**Remark 4.2:** The efficiencies of the estimators discussed in this paper can be compared for fixed cost, following the procedure given in Sukhatme et. al. (1984) and Gupta et. al. (1992-93).

**5. Empirical Study**

To illustrate the performance of various estimators of population correlation coefficient, we consider the data given in Murthy [1967, P.226]. The variates are:
y=output, x=Number of Workers, z =Fixed Capital
N=80, n=10, $n_1$ =25,
$\bar{X} = 283.875$, $\bar{Y} = 5182.638$, $\bar{Z} = 1126$, $C_x = 0.9430$, $C_y = 0.3520$, $C_z = 0.7460$,
$\delta_{003} = 1.030$, $\delta_{004} = 2.8664$, $\delta_{021} = 1.1859$, $\delta_{022} = 3.1522$, $\delta_{030} = 1.295$, $\delta_{040} = 3.65$,
$\delta_{102} = 0.7491$, $\delta_{120} = 0.9145$, $\delta_{111} = 0.8234$, $\delta_{130} = 2.8525$,
$\delta_{112} = 2.5454$, $\delta_{210} = 0.5475$, $\delta_{220} = 2.3377$, $\delta_{201} = 0.4546$, $\delta_{202} = 2.2208$, $\delta_{300} = 0.1301$,
$\delta_{400} = 2.2667$, $\rho_{yx} = 0.9136$, $\rho_{xz} = 0.9859$, $\rho_{yz} = 0.9413$.

The percent relative efficiencies (PREs) of $\hat{\rho}_{1d}, \hat{\rho}_{hd}, \hat{\rho}_{td}$ with respect to conventional estimator r have been computed and compiled in Table 5.1.

Table 5.1: The PRE's of different estimators of $\rho_{yx}$

| Estimator | r | $\hat{\rho}_{hd}$ | $\hat{\rho}_{td}$ (or $\hat{\rho}^*_{td}$) |
|---|---|---|---|
| PRE(.,r) | 100 | 129.147 | 305.441 |

Table 5.1 clearly shows that the proposed estimator $\hat{\rho}_{td}$ (or $\hat{\rho}^*_{td}$) is more efficient than r and $\hat{\rho}_{hd}$.